\documentclass{ifacconf}

\usepackage{natbib}        % required for bibliography

\usepackage{caption}
\usepackage{subcaption}
\usepackage{multirow}

\usepackage{algorithm}
\usepackage{algorithmic}

\usepackage{custom}
\usepackage{ifacthm}

\newcommand{\tr}{\to}
\newcommand{\trs}[1][\sigma]{\tr_{#1}}

\newcommand{\nqu}{n_{q,u}}
\newcommand{\nqx}{n_{q,x}}

\newcommand{\Ps}{\mathcal{P}}
\newcommand{\U}{\mathcal{U}}

\newcommand{\sys}{S}
\newcommand{\new}[1]{#1}

\newcommand{\dtahcs}{HCS}
\newcommand{\dtahas}{HAS}

\usepackage{siunitx}

\usepackage{tkz-graph}
\usepackage{tikz}
\usetikzlibrary{arrows,matrix,decorations.pathreplacing,positioning,chains,fit,shapes,calc} %Voir 1.3.8
\usetikzlibrary{calc,patterns,decorations.pathmorphing,decorations.markings}
\tikzstyle{spring}=[thick,decorate,decoration={zigzag,pre length=0.1cm,post length=0.1cm,segment length=6}]
\tikzstyle{damper}=[thick,decoration={markings,
  mark connection node=dmp,
  mark=at position 0.5 with
  {
    \node (dmp) [thick,inner sep=0pt,transform shape,rotate=-90,minimum width=10pt,minimum height=2pt,draw=none] {};
    \draw [thick] ($(dmp.north east)+(2pt,0)$) -- (dmp.south east) -- (dmp.south west) -- ($(dmp.north west)+(2pt,0)$);
    \draw [thick] ($(dmp.north)+(0,-4pt)$) -- ($(dmp.north)+(0,4pt)$);
  }
}, decorate]

\begin{document}

\begin{frontmatter}

%\clubpenalty=10000
%\widowpenalty = 10000

\title{Computing controlled invariant sets for hybrid systems with applications to model-predictive control}
%\title{On controlled invariant sets for Hybrid Systems : a Sum-of-Squares approach}
%\subtitle{[Extended Abstract]}

\author[First]{Beno\^{i}t Legat}
\author[Second]{Paulo Tabuada}
\author[Third]{Rapha\"{e}l M. Jungers}

\newcommand{\email}[1]{\emph{\texttt{#1}}}

\address[First]{ICTEAM,
  Universit\'e catholique de Louvain,
  4 Av. G. Lema\^itre,
  1348 Louvain-la-Neuve, Belgium
  (e-mail: \email{benoit.legat@student.uclouvain.be})}

\address[Second]{Department of \new{Electrical and Computer} Engineering,
  UCLA,
  (e-mail: \email{tabuada@ee.ucla.edu})}

\address[Third]{ICTEAM,
  Universit\'e catholique de Louvain,
  4 Av. G. Lema\^itre,
  1348 Louvain-la-Neuve, Belgium
  (e-mail: \email{raphael.jungers@uclouvain.be})}

\begin{abstract}
  \,\,\,\,In this paper, we develop a method for computing controlled invariant sets using Semidefinite Programming.
  We apply our method to the controller design problem for switching affine systems with polytopic safe sets.
  The task is reduced to a semidefinite programming problem by enforcing \new{an} invariance relation in the dual space of the geometric problem.
  The paper ends with an application to safety critical model predictive control.
  %We introduce the method for computing controlled invariant sets of a discrete hybrid control affine systems with autonomous switching and polytopic invariant sets.
  %The problem can be solved either using semidefinite, second order cone or linear programming thanks to the DSOS and SDSOS relaxations.
\end{abstract}

\begin{keyword}
  Controller Synthesis; Set Invariance; LMIs; Scalable Methods.
\end{keyword}

\end{frontmatter}

\section{Introduction}
The problem of computing a controlled invariant set is a paradigmatic challenge in the broad field of Hybrid Systems control.
Indeed, it is for instance crucial in safety-critical applications,
such as the control of a platoon of vehicles or air traffic management; see \cite{tomlin1998conflict},
where firm guarantees are needed on our ability to maintain the state in a safe region
(e.g.\new{,} with a certain minimal distance between vehicles).
In other situations, the dynamical system might be too complicated to analyze exactly in every point of the state space,
but yet it can be possible to confine the state within a guaranteed set.
Such situations occur frequently in hybrid, embedded, event-triggered systems, because of the complexity of the dynamics.

%si tu trouves une ou deux références pour les deux applications que je mentionne, qui parlent d'invariant sets (ou si t'as d'autres appli), c'est mieux.  Vu l'urgence je n'ai pas eu le temps de chercher, mais bon on a deux semaines pour faire ça hein...

A set is \emph{controlled invariant} (sometimes also referred to as \emph{viable})
if, any trajectory whose initial point is in the set can be kept inside it by means of a proper control action.
Given a system with constraint specifications on the states and/or input, the controlled invariant set can be used
to determine initial states \new{such that trajectories with these initial conditions are guaranteed to} meet the specifications.
Moreover, \new{in some situations,} a state feedback control law can be derived from the knowledge of the controlled invariant set; see \cite{blanchini1999survey} for a survey.

The computation of invariant sets is usually \new{achieved using either} polyhedral computations or semidefinite programming.
Polyhedral computations are typically restricted to affine constraint specifications but it has been \new{recently shown} that it can also be applied to algebraic constraints; see~\cite{athanasopoulos2016computing}.
If the system contains a control input, the computational complexity of the problem becomes even more challenging.
Indeed, this requires (see e.g., the procedure p.~201 in \cite{blanchini2015set}) the computation of projections of polytopes when using polyhedral computations
and semidefinite programming techniques are not directly applicable.

Methods based on polyhedral computations for hybrid control systems have been developped in \cite{rungger2013specification, smith2016interdependence, rungger2017computing}.
Unfortunately, the problem of polyhedral projection is well known to severely suffer from the curse of dimensionality, see \cite{avis1995good}, and the additional complexity of the discrete dynamics in hybrid systems makes the problem even less scalable for these systems.

The semidefinite programming approach sacrifices exactness of the solution for the sake of algorithmic tractability.
In the case of an uncontrolled system $x_{k+1} = Ax_k$, it consists in
%The semidefinite programming approach for uncontrolled system $x_{k+1} = Ax_k$ consists in
searching for an ellipsoidal set
\[ \mathcal{E}_P = \{\, x \in \R^n \mid x^\Tr P x \leq 1 \,\} \]
such that if $x^\Tr P x \leq 1$ then $x^\Tr A^\Tr P A x \leq 1$.
Indeed, one can verify that it implies invariance of the set $\mathcal{E}_P$.
The S-procedure allows to formulate the search of $P$ as a semidefine program; see \cite{polik2007survey} for a survey on the S-procedure.

With the presence of the control $u$ in the system $x_{k+1} = Ax_k + Bu_k$, the condition becomes:
\[
  x^\Tr P x \leq 1 \Rightarrow \exists u, (A x + Bu)^\Tr P (A x + Bu) \leq 1.
\]
The control term $u$, or more precisely the existential quantifier $\exists$ prevents the S-procedure to be directly applied.

\new{\cite{kurzhanski2005verification} show how to compute an over- and under-approximation of the reachable sets of a hybrid control system.
While they approximate \emph{reachable sets} and do not compute \emph{controlled invariant sets}, their approach bears similarities with the method presented in this paper.
However, their technique does not rely on semidefinite programming as they propagate ellipsoidal sets and do not need to enforce any invariance property.}

%For each ellipsoid, $\mathcal{E}_P$, they store the matrix $P^{-1}$ of the polar set $\mathcal{E}_P^\circ = \mathcal{E}_{P^{-1}}$ and use \lemref{podu} and \eqref{eq:ATr} to compute the image of $\mathcal{E}_P$ under a linear map $A$ as the polar of $\mathcal{E}_{AP^{-1}A^\Tr}$.
%Their method is generalized to ellipsoids not centered at the origin.

In \cite{korda2014convex}, a semidefinite programming method is proposed for the computation of an outer approximation of the maximal controlled invariant sets.
While the set computed with this method can be a good approximation of the maximal controlled invariant set, it is an outer approximation and is not controlled invariant unless the approximation is exact.

In this paper, we give a general method that circumvents this issue.
A key ingredient in our technique is that we work in the dual space of the geometric problem.
We detail the application of the method to two classes of hybrid systems:
Discrete-Time Affine Hybrid Control System (\new{\dtahcs{} for short}) and Discrete-Time Affine Hybrid Algebraic System (\dtahas{} for short).
\dtahas{} are not control systems but the computation of invariant sets for such systems presents
the same features than for \dtahcs{}.
As a matter of fact, we show how to reduce the computation of controlled invariant sets \new{for} \dtahcs{} to the computation of invariant sets \new{for} \dtahas{}.

%Invariant sets are usually computed either as a polytope using polyhedral computation or as the sublevel set of a lyapunov function using

%A controlled invariant set for a time invariant control system is a set such that there exists a state feedback control law such that any trajectory starting inside the set remains in the set.
%The computation of such set has been subject to extensive research over the past decades.

% TODO check \cite{daryin2013parallel}

In this paper we break the problem into four subproblems, which we solve separately.
%The paper is a succession of small propositions, which, put together, allow us to tackle this important problem.
In \secref{constraints}, we show how to reduce the computation of controlled invariant sets of a \dtahcs{} with \emph{constrained} input to controlled invariant sets of a \dtahcs{} with \emph{unconstrained} input.
Then in \secref{dtahas}, we give the reduction of the computation of controlled invariant sets of a \dtahcs{} with unconstrained input to invariant sets of a \dtahas{}.
In \secref{duality}, we detail the relation between the algebraic invariance condition of \dtahas{} on a convex set and its polar set
and we discuss how to lift the state space to handle non-homogeneity.
In \secref{ellipsoids}, we show that using the results of \secref{duality}, the invariance of ellipsoids for a \dtahas{} can be formulated as a semidefinite program.
%We then describe how to use this relation to
%we show how to reduce the computation of controlled invariant sets of a \dtahcs{} with uncconstrained input to invariant sets of a \dtahas{}.

\new{We end the paper with an application of the ellipsoidal controlled invariant sets to safety critical model predictive control.
We show that precomputing such sets allows to guarantee safety of the model predictive controller and thus to alleviate expensive long-horizon computations thereby removing the need for long horizon.}

\section{Controlled Invariant Set}
In this section, we define \dtahcs{} and \dtahas{} and give the invariance conditions for these two classes of hybrid systems.
We detail the relation between controlled invariant sets of \dtahcs{} and invariant sets of \dtahas{}.

\subsection{Discrete-Time Affine Hybrid Control System}
\label{sec:dtahcs}

We will consider the following definition of Discrete-Time Affine Hybrid Control System.

\begin{mydef}
  \label{def:dtahcs}
  A \emph{Discrete-Time Affine Hybrid Control System (\dtahcs{})} is a system $\sys = (T, (A_\sigma, B_\sigma, c_\sigma)_{\sigma \in \Sigma},\\
  (\Ps_q, U_q)_{q \in \Nodes})$ where $T = (\Nodes, \Sigma, \tr)$ \new{and $\tr \subseteq \Nodes \times \Sigma \times \Nodes$}.
  \new{A trajectory is a sequence $\{(x_k,u_k,\sigma_k)\}_{k \in \mathbb{N}}$} satisfying for all $k \in \mathbb{N}$:
  \begin{align*}
    x_{k+1} & = A_{\sigma_{\new{k}}} x_k + B_{\sigma_{\new{k}}} u_k + c_{\sigma_{\new{k}}},\\
    x_k \in \Ps_{q_k}, u_k &\in \U_{q_k}, q_k \trs[\sigma_{\new{k}}] q_{k+1}.
  \end{align*}
  %The system is coupled with a mixed autonomous and controlled switching represented by $\sw = ((\As_q)_{q \in \Nodes}, (\Cs_\as)_{\as})$.
  %The system is coupled with an autonomous switching represented by $\sw = ((\As_q)_{q \in \Nodes}, (\Cs_\as)_{\as})$.
  %Given a current \node{} $q$, a intermediate switching $\as \in \As_q$ is chosen autonomously;
  %then the switching $q \trs q' \in \Cs_{\as}$ is controlled.
\end{mydef}

Given a \node{} $q \in \Nodes{}$, we denote
the set of allowed switching signals as $\Sigma_q$,
the state dimension as $n_{q,x}$ and the input dimension as $n_{q,u}$.

\begin{figure}[!h]
  \centering
  \begin{tikzpicture}
    \draw[thick] (0, 0) -- (2, 0) -- (2, 1) -- (0.5, 1) -- (0, 0.5) -- (0, 0);
    \draw [spring] (2, 0.8) -- (3, 0.8);
    \draw [damper] (2, 0.2) -- (3, 0.2);
    \draw[thick] (3, 0) rectangle (5, 1);
    \draw [spring] (5, 0.8) -- (6, 0.8);
    \draw [damper] (5, 0.2) -- (6, 0.2);
    \draw[thick] (6, 0) rectangle (8, 1);
    \draw[thick] (0.5, 0.1) circle (0.3);
    \fill[color=lightgray] (0.5, 0.1) circle (0.3);
    \draw[thick] (1.5, 0.1) circle (0.3);
    \fill[color=lightgray] (1.5, 0.1) circle (0.3);
    \draw[thick] (3.5, 0.1) circle (0.3);
    \fill[color=lightgray] (3.5, 0.1) circle (0.3);
    \draw[thick] (4.5, 0.1) circle (0.3);
    \fill[color=lightgray] (4.5, 0.1) circle (0.3);
    \draw[thick] (6.5, 0.1) circle (0.3);
    \fill[color=lightgray] (6.5, 0.1) circle (0.3);
    \draw[thick] (7.5, 0.1) circle (0.3);
    \fill[color=lightgray] (7.5, 0.1) circle (0.3);
    \node at (1, 0.5) {$v_0$};
    \node at (4, 0.5) {$v_1$};
    \node at (7, 0.5) {$v_2$};
    \node at (2.5, 1.1) {$d_1$};
    \node at (5.5, 1.1) {$d_2$};
  \end{tikzpicture}
  \caption{Illustration for \exemref{cruise1} with two trailers.}
  \label{fig:cruise}
\end{figure}
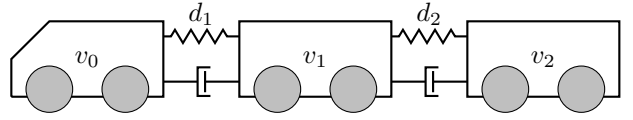

We illustrate this definition with the cruise control example of \cite{rungger2013specification}.
\begin{myexem}
  \label{exem:cruise1}
  We consider a truck with $M$ trailers as represented by \figref{cruise}.
  There is a truck with mass $m_0$ and speed $v_0$ followed by multiple trailers with mass $m$ each.
  The speed of the $i$th trailer is denoted $v_i$.
  There is a spring with stiffness $k_d$ and elongation $d_1$ (resp. $d_i$) and a damper with coefficient $k_s$
  between the truck and the first trailer (resp. the $(i-1)$th trailer and the $i$th trailer).
  The scalar input $u$ controls the speed $v_0$ of the truck by creating a force $m_0u$.
  The dynamics of the system is given by the following equations:
  \begin{align}
    \notag
    \dot{v}_0 & = \frac{k_d}{m_0}(v_1 - v_0) - \frac{k_s}{m_0} d_1 + u\\
    \notag
    \dot{v}_i & = \frac{k_d}{m}(v_{i-1} - 2v_i + v_{i+1}) + \frac{k_s}{m} (d_i - d_{i+1}) & 1 \leq i < M\\
    \label{eq:truckdyn}
    \dot{v}_M & = \frac{k_d}{m}(v_{M-1} - v_M) + \frac{k_s}{m} d_M\\
    \notag
    \dot{d}_i & = v_{i-1} - v_i & 1 \leq i \leq M.
  \end{align}
  %where the entries of the state $x = (d, v_1, v_0) \in \R^3$ correspond to the distance $d$ between the truck and the trailer, the speed of the trailer $v_1$ and the speed of the truck $v_0$.
  %The constraints on the state space are as follows: the state should...
  The spring elongation should always remain between $-\SI{0.5}{\meter}$ and $\SI{0.5}{\meter}$ and the speeds of the truck and trailers should remain between $\SI{5}{\meter\per\second}$ and $\SI{35}{\meter\per\second}$.
  Moreover, there are three speed limits $\bar{v}_a = \SI{15.6}{\meter\per\second}$, $\bar{v}_b = \SI{24.5}{\meter\per\second}$, $\bar{v}_c = \SI{29.5}{\meter\per\second}$ and
  whenever the truck is informed of a new speed limit, it has \SI{0.8}{\second} to decrease $v_i$ ($0 \leq i \leq M$) below the speed limit.

  We sample time with a period of \SI{0.4}{\second} and define an initial \node{} $q_{d0}$ and 6 \nodes{} $q_{ij}$ where $i \in \{a, b, c\}$ is the current speed limitation and $j \in \{0, 1\}$ is the number of sampling times left to satisfy the limit.
  The transitions are $q_{ij} \trs[\sigma] q_{\sigma 1}$ for each $i \in \{a, b, c, d\}$ and $\sigma \in \{a, b, c, d\} \setminus \{i\}$.
  The symbol $a$ (resp. $b$, $c$) represents that the truck sees a new speed limitation $\bar{v}_a$ (resp. $\bar{v}_b$, $\bar{v}_c$) and $d$ represents that it does not see any new speed limitation.
  We suppose for simplicity that it is not possible to see a new speed limitation $\bar{v}_\sigma$ from a \node{} $q_{\sigma j}$.
  %\Cs_i & = \{\, q_{kj} \trs[1] q_{i1} \mid j = 0, 1 \,\} \quad i = a, b, c\\
%  The autonomous switching sets are $\As_{d0} = \{a, b, c, d\}$ and for $i = a, b, c$, $\As_{q_{ij}} = \{a, b, c, d\} \setminus \{i\}$ where $a$ (resp. $b$, $c$) represents that the truck sees a new speed limitation $\bar{v}_a$ (resp. $\bar{v}_b$, $\bar{v}_c$) and $d$ represents that it does not see any new speed limitation.
%  We exclude $i$ from $\As_{q_{ij}}$ for simplicity since it is equivalent to $d$.
  \new{The possible transitions are} represented in \figref{cruise1}.

  \begin{figure}[!ht]
    \centering
    \begin{tikzpicture}[x=40,y=18]
      \Vertex[L={$q_{d0}$}]{0}
      \NOEA[unit=2,L={$q_{a1}$}](0){1}
      \EA[unit=2,L={$q_{a0}$}](1){2}
      %\EA[unit=2,L={b1}](0){3}
      %\EA[unit=2,L={b0}](3){4}
      \SOEA[unit=2,L={$q_{c1}$}](0){5}
      \EA[unit=2,L={$q_{c0}$}](5){6}
%     \SOEA[unit=2](2){3}
%     \EA[unit=2](3){4}
      \tikzset{EdgeStyle/.style = {->}}
      \tikzset{LoopStyle/.style = {->}}
      \Edge[label=$a$](0)(1)
      \Edge[label=$d$](1)(2)
      %\Edge[label=$b$](0)(3)
      %\Edge[label=$d$](3)(4)
      \Edge[label=$c$](0)(5)
      \Edge[label=$d$](5)(6)
      %\tikzset{EdgeStyle/.append style = {bend right=18}}
      %\Edge[label=$a$](3)(1)
      %\Edge[label=$b$](1)(3)
      %\Edge[label=$c$](3)(5)
      %\Edge[label=$b$](5)(3)
      \tikzset{EdgeStyle/.append style = {bend right=20}}
      \Edge[label=$c$](2)(5)
      \tikzset{EdgeStyle/.append style = {bend right=-20}}
      \Edge[label=$a$](6)(1)
      \tikzset{EdgeStyle/.append style = {bend right=10}}
      \Edge[label=$c$](1)(5)
      \Edge[label=$a$](5)(1)
%     \Edge[label=3](1)(2)
%     \Edge[label=2](2)(1)
%     \Edge[label=1](2)(3)
%     \Edge[label=2](3)(1)
%     \Edge[label=3](3)(2)
%     \Edge[label=4](3)(4)
      %\draw[thick,->] (0) to [out=65,in=115,looseness=10] node [midway, fill=white] {$d$} (0);
      \draw[thick,->] (2) to [out=-115,in=-65,looseness=9] node [midway, fill=white] {$d$} (2);
      \draw[thick,->] (6) to [out=65,in=115,looseness=9] node [midway, fill=white] {$d$} (6);
%     \Edge[label=1](4)(3)
    \end{tikzpicture}
     \caption{Transitions and switchings between the \nodes{} for \exemref{cruise1}.
     Nodes $q_{b1}$ and $q_{b0}$ are not shown for clarity.
     }
     \label{fig:cruise1}
  \end{figure}

%  The controlled switching are
%  \begin{align*}
%    \Cs_i & = \{\, q_{kj} \trs[1] q_{i1} \mid j = 0, 1 \,\} \quad i = a, b, c\\
%    \Cs_d & = \{\, q_{ij} \trs[1] q_{i(j-1)} \mid i = a, b, c, j = 1, 2 \,\}.
%  \end{align*}
%  Note that for each $q_{ij}$, $\alpha \in \As_{q_{ij}}$, there is only one transition
%  in $\Cs_\alpha$ that starts at $q_{ij}$ so the switching can be considered purely autonomous.

  %The reset map 0 is defined by $A_0 = I$, the identity matrix, $B_0 = 0$ and $c_0 = 0$.
  The reset maps $(A_\sigma, B_\sigma, c_\sigma)$ are simply the integration of the dynamical system \eqref{eq:truckdyn} over \SI{0.4}{\second} with a zero-order hold input extrapolation.

  Let
  \begin{align*}
    P_0 & = \{\, (d, v) \in \R^{2M+1} \mid -0.5 \leq d \leq 0.5, 5 \leq v \leq 35 \,\},\\
    P_i & = \{\, (d, v) \in \R^{2M+1} \mid v \leq \bar{v}_i \,\}, \quad i = a, b, c,
  \end{align*}
  where $d = (d_1, \ldots, d_M)$, $v = (v_0, \ldots, v_M)$ and inequalities in the two equations above are entrywise.
  The safe sets are $\Ps_{q_{d0}} = P_0$ and for $i = a, b, c$, $P_{q_{ij}} = P_0$ if $j > 0$ and $P_{q_{i0}} = P_0 \cap P_i$.
  The input set is $\U_{ij} = \{\, u \in \R \mid -4 \leq u \leq 4 \,\}$ for each \node{} $q_{ij}$.
\end{myexem}

\begin{mydef}[Controlled invariant sets for a \dtahcs{}]
  \label{def:cis}
  $ $\\
  Consider a \dtahcs{} $\sys$.
  We say that sets $\Csetvar = (\Csetvar_q)_{q \in \Nodes}$ are \emph{controlled invariant} for $\sys$ if
  $\Csetvar_q \subseteq \Ps_q$ for each $q \in \Nodes$ and
  %$\forall q \in \Nodes, x \in \Csetvar_q, \as \in \As_q$, $\exists q \trs q' \in \Cs_{\as}, u \in \U_v$ such that
  $\forall x \in \Csetvar_q, q \trs[\sigma] q'$, $\exists u \in \U_q$ such that
  \[ A_\sigma x + B_\sigma u + c_\sigma \in \Csetvar_{q'}. \]
\end{mydef}

\begin{myrem}
  \label{rem:autonomous}
  It is important to distinguish two types of switching:
  \emph{autonomous switching} and \emph{controlled switching};
  see details in \cite[Section~1.1.3]{liberzon2012switching}.
  \defref{cis} is the definition of controlled invariance for autonomous systems and in this paper we only consider systems that switch autonomously.
  With controlled switching, ``$\forall q \trs[\sigma] q'$'' is replaced by ``$\exists q \trs[\sigma] q'$'' in \defref{cis}.
\end{myrem}

\subsection{Handling controller constraints}
\label{sec:constraints}
We say that the input of a \dtahcs{} is \emph{unconstrained} if $\U_q = \R^{\nqu}$ for all $q \in \Nodes$, otherwise we say that the input is \emph{constrained}.
The computation of controlled invariant sets for a \dtahcs{} with constrained input can be reduced to the computation of invariant sets for a \dtahcs{} with unconstrained input
as shown by the following lemma.

\begin{mylem}
  \label{lem:liftu}
  The sets $\Csetvar = (\Csetvar_q)_{q \in \Nodes}$ are \emph{controlled invariant} for $\sys = (T, (A_\sigma, B_\sigma, c_\sigma)_{\sigma \in \Sigma}, (\Ps_q, U_q)_{q \in \Nodes})$ if and only if
  their exist controlled invariant sets $\Csetvar' = (\Csetvar_q')_{q \in \Nodes'}$ such that $\Csetvar'_q = \Csetvar_q$ $\forall q \in \Nodes$ for the system
  $\sys' = (T', (A_\sigma, B_\sigma, c_\sigma)_{\sigma \in \Sigma'}, (\Ps'_q, \U'_q)_{q \in \Nodes'})$
  where the new transitions $T' = (\Nodes', \Sigma', \tr')$ are obtained as follows:
  %For each \node{} $q$ and autonomous switching $\as \in \As_q$, we create a \node{} $q^\as$ and transform each transition $q \trs q' \in \Cs_\as$
  %into $q \trs[q^0] q^\as$ and $q^\as \trs[\sigma'] q'$.
  For each transition $q \trs \new{r}$ \new{in $T$}, we create a \node{} $q^\sigma$
  and the transitions $q \trs[q^0]' q^\sigma$ and $q^\sigma \trs[\sigma']' \new{r}$ \new{in $T'$}.

  The new safe and input sets are
  \begin{align*}
    \Ps'_q & = \Ps_q & \U'_q & = \R^{\nqu}\\
    \Ps'_{q^\sigma} & =  \Ps_q \times \U_q & \U'_{q^\sigma} & = \R^0
  \end{align*}
  and the new reset maps are
  \begin{align*}
    A_{q^0} & =
    \begin{bmatrix}
      I\\
      0
    \end{bmatrix} &
    B_{q^0} & =
    \begin{bmatrix}
      0\\
      I
    \end{bmatrix} &
    c_{q^0} & = 0\\
    A_{\sigma'} & =
    \begin{bmatrix}
      A_\sigma & B_\sigma
    \end{bmatrix} &
    %B_{\sigma'} & =
    %\begin{bmatrix}
    %  \, \, \,
    %\end{bmatrix} &
    &&
    c_{\sigma'} & = c_\sigma
  \end{align*}
  \new{and $B_{\sigma'}$ is the unique map sending $0 \in \R^0$ to $0 \in \mathbb{R}^{n_r}$.}

  \begin{proof}
    Consider controlled invariant sets $\Csetvar'$ for $\sys'$ and let $\Csetvar = (\Csetvar'_q)_{q \in \Nodes}$.
    %Given $q \in \Nodes, x \in \Csetvar_q, \as \in \As_q$,
    %the controlled invariance of $\Csetvar'$ ensures that there exists $u, q'$ such that $(x, u) \in \Csetvar'_{q^\as} \subseteq \Ps_q \times \U_q$, $q \trs q' \in \Cs_\as$ and $A_\sigma x + B_\sigma u + c_\sigma \in \Csetvar'_{q'} = \Csetvar_{q'}$. Hence $\Csetvar$ is controlled invariant for $\sys$.
    Given $x \in \Csetvar_q$ and $q \trs r$,
    the controlled invariance of $\Csetvar'$ ensures that there exists $u$ such that $(x, u) \in \Csetvar'_{q^\sigma} \subseteq \Ps_q \times \U_q$ and $A_\sigma x + B_\sigma u + c_\sigma \in \Csetvar'_{r} = \Csetvar_{r}$. Hence $\Csetvar$ is controlled invariant for $\sys$.

    Consider now controlled invariant sets $\Csetvar$ for $\sys$ and let $\Csetvar' = (\Csetvar'_q)_{q \in \Nodes'}$ where $\Csetvar'_q = \Csetvar_q$ for each $q \in \Nodes$.
    %For each $q \in \Nodes, x \in \Csetvar'_q = \Csetvar_q, \as \in \As_q$, the controlled invariance of $\Csetvar$ ensures that there exists $u \in \U_q$, $q \trs q' \in \Cs_\as$ such that $A_\sigma x + B_\sigma u + c_\sigma \in \Csetvar_{q'} = C'_{q'}$, setting $\Csetvar'_{q^\as}$ to be the union of these pairs $(x, u)$ makes $\Csetvar'$ controlled invariant for $\sys'$.
    Given $q \trs r$, for each $x \in \Csetvar'_q = \Csetvar_q$ the controlled invariance of $\Csetvar$ ensures that there exists $u \in \U_q$ such that $A_\sigma x + B_\sigma u + c_\sigma \in \Csetvar_{r} = \Csetvar'_{r}$, setting $\Csetvar'_{q^\sigma}$ to be the union of these pairs $(x, u)$ makes $\Csetvar'$ controlled invariant for $\sys'$.
  \end{proof}
\end{mylem}

\begin{myrem}
  \label{rem:liftu}
  %If for a given $q$, $\As_q$ is a singleton $\{\as\}$, we can merge $q$ and $q^\as$ into one state hence have $\Ps'_q = \Ps_q \times \U_q$.
  If for a given $q$, $\Sigma_q$ is a singleton $\{\sigma\}$, we can merge $q$ and $q^\sigma$ into one state hence have $\Ps'_q = \Ps_q \times \U_q$.
  In that case, $\Csetvar_q$ will be the projection of $\Csetvar'_q$ in its state space.
  Even if $\Sigma_q$ is not a singleton, we can pick a single $\sigma \in \Sigma_q$ and merge $q$ and $q^\sigma$ into one state
  and use the reset map
  \begin{align*}
    A_{q^0} & =
    \begin{bmatrix}
      I & 0\\
      0 & 0
    \end{bmatrix} &
    B_{q^0} & =
    \begin{bmatrix}
      0\\
      I
    \end{bmatrix} &
    c_{q^0} & = 0\\
  \end{align*}
  so that switchings $\sigma' \in \Sigma_q \setminus \{\sigma\}$ ignore the part of the state of $q$ that corresponds to the input to be used for $\sigma$.

  %The reason we need to create a new state for each $\as \in \As_q$
  %When applying \lemref{liftu}, we need to create a new state for each $q \in \Nodes$ and each $\as \in \As_q$ because
\end{myrem}

\begin{myexem}
  \label{exem:cruise2}
  We represent on \figref{cruise2} the application of the transformation described in \lemref{liftu} to the system of \exemref{cruise1}.
  We can use \remref{liftu} to avoid creating $q^d$ for each $q$.
  Moreover, since $(A_\sigma, B_\sigma, c_\sigma)$ does not depend on $\sigma$, we can merge all the \nodes{} $q^a$ (resp. $q^b$, $q^c$) together
  into a common state that we name $q_{a2}$ (resp. $q_{b2}$, $q_{c2}$).
  %This transformation is represented in \figref{cruise2}.
  %Consider that we want to transform the system of \exemref{cruise1} into a system with unconstrained input.

  %This is represented in \figref{cruise2}.
  \begin{figure}[!ht]
    \centering
    \begin{tikzpicture}[x=40,y=18]%[x=60,y=30]
      \Vertex[L={$q_{d0}$}]{0}
      \NO[unit=2,empty=true](0){13}
      \EA[unit=1,L={$q_{a2}$}](13){12}
      \EA[unit=2,L={$q_{a1}$}](12){11}
      \EA[unit=2,L={$q_{a0}$}](11){10}
      %\EA[unit=2,L={b1}](0){3}
      %\EA[unit=2,L={b0}](3){4}
      \SO[unit=2,empty=true](0){33}
      \EA[unit=1,L={$q_{c2}$}](33){32}
      \EA[unit=2,L={$q_{c1}$}](32){31}
      \EA[unit=2,L={$q_{c0}$}](31){30}
%     \SOEA[unit=2](2){3}
%     \EA[unit=2](3){4}
      \tikzset{EdgeStyle/.style = {->}}
      \tikzset{LoopStyle/.style = {->}}
      \Edge[label=$a$](0)(12)
      \Edge[label=$d$](12)(11)
      \Edge[label=$d$](11)(10)
      %\Edge[label=$b$](0)(3)
      %\Edge[label=$d$](3)(4)
      \Edge[label=$c$](0)(32)
      \Edge[label=$d$](32)(31)
      \Edge[label=$d$](31)(30)
      %\tikzset{EdgeStyle/.append style = {bend right=18}}
      %\Edge[label=$a$](3)(1)
      %\Edge[label=$b$](1)(3)
      %\Edge[label=$c$](3)(5)
      %\Edge[label=$b$](5)(3)
      \tikzset{EdgeStyle/.append style = {bend right=6}}
      \Edge[label=$c$](10)(32)
      \tikzset{EdgeStyle/.append style = {bend right=-6}}
      \Edge[label=$a$](30)(12)
      \tikzset{EdgeStyle/.append style = {bend right=20}}
      \Edge[label=$c$](11)(32)
      \tikzset{EdgeStyle/.append style = {bend right=-20}}
      \Edge[label=$a$](31)(12)
      \tikzset{EdgeStyle/.append style = {bend right=10}}
      \Edge[label=$c$](12)(32)
      \Edge[label=$a$](32)(12)
%     \Edge[label=3](1)(2)
%     \Edge[label=2](2)(1)
%     \Edge[label=1](2)(3)
%     \Edge[label=2](3)(1)
%     \Edge[label=3](3)(2)
%     \Edge[label=4](3)(4)
      %\draw[thick,->] (0) to [out=65,in=115,looseness=10] node [midway, fill=white] {$d$} (0);
      \draw[thick,->] (10) to [out=-115,in=-65,looseness=7] node [midway, fill=white] {$d$} (10);
      \draw[thick,->] (30) to [out=65,in=115,looseness=7] node [midway, fill=white] {$d$} (30);
%     \Edge[label=1](4)(3)
    \end{tikzpicture}
     \caption{Transitions and switchings between the \nodes{} for \exemref{cruise2}.
     \Nnodes{} $q_{b2}$, $q_{b1}$ and $q_{b0}$ are not shown for clarity.
     }
     \label{fig:cruise2}
  \end{figure}
\end{myexem}

\subsection{Discrete-Time Affine Hybrid Algebraic System}
\label{sec:dtahas}

\begin{mydef}
  \label{def:dtahas}
    A \emph{Discrete-Time Affine Hybrid Algebraic System (\dtahas{})} is a system $\sys = (T, (A_\sigma, E_\sigma, c_\sigma)_{\sigma \in \Sigma},\\
    (\Ps_q)_{q \in \Nodes})$ where $T = (\Nodes, \Sigma, \tr)$ \new{and $\tr \subseteq \Nodes \times \Sigma \times \Nodes$}.
  \new{A trajectory is a sequence $\{(x_k,\sigma_k)\}_{k \in \mathbb{N}}$} satisfying for all $k \in \mathbb{N}$:
  \begin{align*}
    E_{\sigma_{\new{k}}} x_{k+1} & = A_{\sigma_{\new{k}}} x_k + c_{\sigma_{\new{k}}},\\
    x_k \in \Ps_{q_k}, u_k &\in \U_{q_k}, q_k \trs[\sigma_{\new{k}}] q_{k+1}.
  \end{align*}
  %The switching $\sw$ is defined exactly as \defref{dtahcs}
\end{mydef}

%We say that a \dtahcs{} or \dtahas{} has \emph{autonomous switching} if $\as$ is a singeton for each $q \in \Nodes$ and $\as \in \As_q$.

\begin{mydef}[Invariant sets for a \dtahas{}]
  \label{def:is}
  Consider a \dtahas{} $\sys$.
  We say that sets $\Csetvar = (\Csetvar_q)_{q \in \Nodes}$ are \emph{invariant} for $\sys$ if
  $\Csetvar_q \subseteq \Ps_q$ for each $q \in \Nodes$ and
  %$\forall q \in \Nodes, x \in \Csetvar_q, \as \in \As_q$, $\exists q \trs q' \in \Cs_{\as}$ such that
  for all $q \trs q'$,
  \begin{equation}
    \label{eq:isdtahas}
    A_\sigma \Csetvar_{q} + c_\sigma \subseteq E_\sigma \Csetvar_{q'}.
  \end{equation}
\end{mydef}

\begin{myrem}
  \label{rem:descriptor}
  \defref{is} \new{can be interpreted as stating that} $\Csetvar$ is invariant if for each transition $q \trs q'$ and $x \in \Csetvar_q$,
  \begin{quote}
    there \emph{exists} $y \in \Csetvar_{q'}$ such that $A_\sigma x + c_\sigma = E_\sigma y$.
  \end{quote}
  A similar definition exists where this last part is replaced by
  \begin{quote}
    for \emph{each} $y$ such that $A_\sigma x + c_\sigma = E_\sigma y$, $y$ must belong to $\Csetvar_{q'}$.
  \end{quote}
  This is not equivalent to \defref{is} if $A_\sigma$ and $E_\sigma$ are not full rank.
  Moreover, computing ellipsoidal invariant sets according to this definition \new{is much easier}: it simply amounts to finding
  positive definite matrices $Q_q$ such that $A_\sigma^\Tr Q_q A_\sigma \preceq E_\sigma^\Tr Q_{q'} E_\sigma$; see \cite{owens1985consistency}.
\end{myrem}

We now show that the computation of controlled invariant sets of a \dtahcs{} can be reduced to
the computation of invariant sets of a \dtahas{}.

\begin{mylem}
  \label{lem:proju}
  The sets $\Csetvar = (\Csetvar_q)_{q \in \Nodes}$ are \emph{controlled invariant} for the \dtahcs{}
  $\sys = (T, (A_\sigma, B_\sigma, c_\sigma)_{\sigma \in \Sigma}, (\Ps_q, \R^{\nqu})_{q \in \Nodes})$ if and only if
  they are invariant sets for the \dtahas{}
  $\sys' = (T, (E_\sigma A_\sigma, E_\sigma, E_\sigma c_\sigma)_{\sigma \in \Sigma}, (\Ps_q)_{q \in \Nodes})$
  where $E_\sigma$ is a projection on $\Image(B_\sigma)^{\perp}$.

  \begin{proof}
    As the input is unconstrained,
    for each $q \trs q'$ and $x \in \Ps_q$,
    there exists $u \in \R^{n_{q,u}}$ such that $A_\sigma x + B_\sigma u + c_\sigma \in \Csetvar_{q'}$ if and only if
    $E_\sigma A_\sigma x + E_\sigma c_\sigma \in E_\sigma \Csetvar_{q'}$.
  \end{proof}
\end{mylem}

% TODO Compare with liberzon2012switched and liberzon2009stability

%Note that if $\Csetvar$ and $\Csetvar'$ are controlled invariant sets of a \dtahcs{} (resp. invariant sets of a \dtahas{}) then their respective \emph{unions} $\Csetvar'' = (\Csetvar_q \cup C'_q)_{q \in \Nodes}$ are also controlled invariant sets (resp. invariant sets) of the system.
%
%\begin{mydef}
%  The \emph{maximal controlled invariant} sets of a \dtahcs{} (resp. \emph{maximal invariant} sets of a \dtahas{}) are the union of all the controlled invariant sets (resp. invariant sets) of the system. \todo[inline]{We do not really use this so it can be removed if we need more space}
%\end{mydef}

%\begin{myprob}
%  Given a specification provided as an automaton $G$, find a controlled invariant set for system \defref{dtahcs}.
%\end{myprob}

\section{Computing controlled invariant sets}

\subsection{Duality correspondence for the invariance condition}
\label{sec:duality}

Given a set $\Csetvar$ and a linear map $A$, we define the following notations:
\begin{align}
  \notag
  A\Csetvar & = \{\, Ax \mid x \in \Csetvar \,\}\\
  \notag
  A^{-1}\Csetvar & = \{\, x \mid Ax \in \Csetvar \,\}\\
  \label{eq:ATr}
  A^{-\Tr}\Csetvar & = \{\, x \mid A^{\Tr}x \in \Csetvar \,\}.
\end{align}
Note that $A$ does not need to be invertible in these definitions.

Invariant sets can be computed numerically as \emph{sublevel sets}\footnote{\new{The $\ell$-sublevel set of a function $f : \mathbb{R}^n \to \mathbb{R}$ is the set $\{\, x \in \R^n \mid f(x) \leq \ell \,\}$.}} of polynomials functions using Sum-of-Squares.
One property of sublevel sets that is usually used can be formulated as follows:
If $\Csetvar$ is the $\ell$-sublevel set of a function $f$ then for any function $g$,
$g^{-1}(\Csetvar)$ is the $\ell$-sublevel set of the function $f \circ g$.
Thanks to this property, computing a set $\Csetvar$ satisfying $A\Csetvar \subseteq \Csetvar$ for some linear map $A$ can be for example achieved by
searching for a set $\Csetvar$ being the 1-sublevel set of a polynomial $p(x)$.
Indeed, %using the S-procedure,
the invariance constraint is equivalent to $\Csetvar \subseteq A^{-1}\Csetvar$ which is equivalent to the following \new{implication} : for all $x$, $p(x) \leq 1 \Rightarrow p(\new{A}x) \leq 1$.
The latter proposition can be translated to a constraint of nonnegativity of a polynomial using the Sum-of-Squares formulation and the S-procedure.

\begin{mylem}[S-procedure]
  \label{lem:sproc}
  Given two symmetric matrices $Q_1, Q_2 \in \R^{n \times n}$,
  the existence of a $\lambda \geq 0$ such that the matrix $\lambda Q_{\new{1}} - Q_{\new{2}}$ is positive semidefinite
  is sufficient for the following proposition to hold:
  \begin{quote}
    for all $x \in \R^n$, $x^\Tr Q_1 x \leq 0 \Rightarrow x^\Tr Q_2 x \leq 0$
  \end{quote}
  Moreover, if there exists $x \in \R^n$ such that $x^\Tr Q_1 x > 0$ then this condition is also necessary.
\end{mylem}
%the existence of a number $\lambda > 0$ such that $\lambda p(x) \geq p(Ax)$.

For \dtahas{}, we have in \eqref{eq:isdtahas} an invariance constraint of the form $A\Csetvar \subseteq E\Csetvar$ and we would like to find an equivalent form with
a pre-image as we had with $\Csetvar \subseteq A^{-1}\Csetvar$.
This can be achieved using the polar of the set $\Csetvar$ thanks to the following lemma.
\begin{mylem}[{\cite[Corollary~16.3.2]{rockafellar2015convex}}]
  \label{lem:podu}
  For any convex set $\Csetvar$ (resp. convex cone $\Kset$) and linear map $A$,
  \begin{align*}
    (A\Csetvar)^{\circ} & = A^{-\Tr} \Csetvar^\circ\\
    (A\Kset)^* & = A^{-\Tr} \Kset^*
  \end{align*}
  where $\Csetvar^\circ$ denotes the polar of the set $\Csetvar$ and $\Kset^*$ denotes the dual of the cone $\Kset$.
\end{mylem}

\lemref{podu} shows that $A\Csetvar \subseteq E\Csetvar$ is equivalent to $A^{-\Tr}\Csetvar^\circ \supseteq E^{-\Tr}\Csetvar^\circ$.
Since the invariant sets of the \dtahas{} may not have the origin in their interior, the polar transformation cannot be readily applied.
We handle this non-homogeneity by taking the conic hull of the lifted sets $\Csetvar \times \{1\}$.
More precisely, we define
\begin{align}
  \tau(\Csetvar) & = \{\, (\lambda x, \lambda) \mid \lambda \geq 0, x \in \Csetvar \,\}\\
  r(A, c) & =
  \begin{bmatrix}
    A & c\\
    0 & 1
  \end{bmatrix}.
\end{align}

It can be verified that for any set $\Csetvar$, vector $c$ and linear map $A$,
\begin{equation}
  \label{eq:taur}
  \tau(A\Csetvar + c) = r(A, c) \tau(\Csetvar).
\end{equation}
Moreover, for any half-space $a^\Tr x \leq \beta$,
\begin{equation}
  \label{eq:hs}
  a^\Tr x \leq \beta, \forall x \in \Csetvar \Leftrightarrow (-a, \beta) \in \tau(\Csetvar)^*.
\end{equation}

\begin{mytheo}
  \label{theo:dtahas}
  Consider a \dtahas{} $\sys$. % that has autonomous switching.
  \new{The} closed convex sets $\Csetvar = (\Csetvar_q)_{q \in \Nodes}$ are \emph{invariant} for $\sys$ if and only if
  $\Csetvar_q \subseteq \Ps_q$ for each $q \in \Nodes$ and
  %$\forall q \in \Nodes, \as \in \As_q$,
  %for the only element $q \trs q' \in \Cs_{\as}$,
  for all $q \trs q'$,
  \begin{equation}
    \label{eq:dtahas}
    r(A_\sigma, c_\sigma)^{-\Tr} \tau(\Csetvar_q)^* \supseteq r(E_\sigma, 0)^{-\Tr} \tau(\Csetvar_{q'})^*.
  \end{equation}
  \begin{proof}
    The invariance constraint of \defref{is}
    \[ A_\sigma \Csetvar_q + c_\sigma \subseteq E_\sigma \Csetvar_{q'} \]
    can be rewritten, using \eqref{eq:taur}, into
    \begin{equation}
      \label{eq:coneinv}
      r(A_\sigma, c_\sigma) \tau(\Csetvar_q) \subseteq r(E_\sigma, 0) \tau(\Csetvar_{q'}).
    \end{equation}
    As the sets $\Csetvar_q$ are closed and convex, so are the cones $\tau(\Csetvar_q)$
    hence $\tau(\Csetvar_q)^{**} = \tau(\Csetvar_q)$.
    Therefore, by \lemref{podu}, \eqref{eq:coneinv} is equivalent to \eqref{eq:dtahas}.
  \end{proof}
\end{mytheo}

\subsection{Computation using ellipsoids}
\label{sec:ellipsoids}

While \theoref{dtahas} holds for any convex sets $(\Csetvar_q)_{q \in \Nodes}$,
restricting our attention to ellipsoidal sets renders the invariance condition \eqref{eq:coneinv} amenable to semidefinite programming.
Using sublevel sets of polynomials of higher degree would also allow us to use semidefinite programming but we only describe the ellipsoidal case for simplicity.
This section details the semidefinite program needed to find these ellipsoidal invariant sets and shows its exactness in \theoref{quad}.
%In this section, we show how to compute ellipsoidal invariant sets for \dtahas{} (hence also controlled invariant sets for a \dtahcs{} thanks to \lemref{liftu} and \lemref{proju}).
%The choice of ellipsoid for the shape of our invariant sets allows us to use semidefinite programming.

%\theoref{dtahas} can readily be used to compute ellipsoidal invariant sets for a \dtahas{}.

We define the following notations for ellipsoids
\begin{align*}
  \Ellc{Q}{c} & = \{\, x \mid (x-c)^\Tr Q (x-c) \leq 1 \,\}\\
  \Ellq{D}{d}{\delta} & = \{\, x \mid x^\Tr D x + 2d^\Tr x + \delta \leq 0 \,\}.
\end{align*}

We denote the set of symmetric matrices of $\R^{n}$ as $\SymK$.
\begin{mylem}
  \label{lem:QD}
  Let $Q, D \in \SymK$, $c, d \in \R^n$, $\delta \in \R$ with $Q \succ 0$.
  We have $\Ellc{Q}{c} = \Ellq{D}{d}{\delta}$ if and only if $D \succ 0$ and there exists $\lambda > 0$ such that
  \begin{align}
    \label{eq:lambda}
    \lambda & = d^\Tr D^{-1} d - \delta\\
    \label{eq:c}
    c & = -D^{-1}d\\
    \label{eq:Q}
    Q & = D/\lambda.
  \end{align}
  \begin{proof}
    Substituting $Q$ and $c$ using \eqref{eq:c} and \eqref{eq:Q} in $(x-c)^\Tr Q (x-c) - 1$ gives $(x^\Tr D x + 2d^\Tr x + d^\Tr D^{-1} d - \lambda) / \lambda$.
    We can conclude the ``if'' part of the proof with \eqref{eq:lambda}.
    We now show the ``only if'' part.

    By \lemref{sproc}, for $\Ellc{Q}{c} = \Ellq{D}{d}{\delta}$ to hold, there must exist $\lambda > 0$ such that
    \[ x^\Tr D x + 2d^\Tr x + \delta = \lambda((x-c)^\Tr Q (x-c) - 1). \]
    This implies that
    \begin{align}
      \label{eq:delta}
      \delta & = \lambda c^\Tr Q c - \lambda\\
      \label{eq:d}
      d & = -\lambda Qc\\
      \label{eq:D}
      D & = \lambda Q.
    \end{align}
    Equations~\eqref{eq:d} and \eqref{eq:D} directly give \eqref{eq:c} and \eqref{eq:Q}.
    It remains to show \eqref{eq:lambda}.
    Equation~\eqref{eq:d} is equivalent to $Q^{-1/2}d = -\lambda Q^{1/2}c$ which implies
    \begin{equation}
      \label{eq:dQd}
      d^\Tr Q^{-1} d = \lambda^2 c^\Tr Q c.
    \end{equation}
    Combining \eqref{eq:dQd} with \eqref{eq:D}, we get $\lambda c^\Tr Q c = d^\Tr D^{-1} d$
    which, combined with \eqref{eq:delta}, gives \eqref{eq:lambda}.
  \end{proof}
\end{mylem}

We use the following corollary to represent the cones $\tau(\Csetvar_q)^*$ as the 0-sublevel set of quadratic forms $p(y) = p(x, z) = x^\Tr D_q x + 2d_q^\Tr xz + \delta_q z^2$.

\begin{mycoro}
  \label{coro:convexcone}
  Let $\Kset = \{\, (x, z) | x^\Tr D x + 2d^\Tr xz + \delta z^2 \leq 0, z \geq 0 \,\}$ be a cone that has a nonempty interior
  and no intersection
  with the hyperplane $\{\, (x, 0) | x \in \R^n \,\}$ except the origin.
  The cone $\Kset$ is convex if and only if $D \succ 0$.
  \begin{proof}
    Let $\Csetvar = \Ellq{D}{d}{\delta}$.
    Since every point of the cone satisfy $z > 0$ except the origin, we have $\tau(\Csetvar) = \Kset$.
    Therefore, $\Kset$ is convex if and only if $\Csetvar$ is convex.
    Since $\Kset$ is nonempty,
    \[ \delta - d^\Tr D d = \min_{x \in \R^n} x^\Tr D x + 2d^\Tr x + \delta < 0. \]
    We conclude with \lemref{QD}.
  \end{proof}
\end{mycoro}

%We formulate the search for the cones $\tau(\Csetvar_q)^*$ of \eqref{eq:dtahas} when the sets $\Ps_q$ are polyhedra as a semidefinite programming
%by searching for a cone $\tau(\Csetvar_q)^*$ being the 0-sublevel set of a quadratic form $p(y) = p(x, z) = x^\Tr D_q x + 2d_q^\Tr xz + \delta z^2$.

In \cororef{convexcone}, we require the cone to have no intersection with an hyperplane (except the origin).
However, the cone $\tau(\Csetvar_q)^*$ has no intersection with the hyperplane $\{\, (x, 0) | x \in \R^n \,\}$ if and only if the origin is contained in $\Csetvar_q$ which may not be the case.
In order to alleviate this, the approach we suggest is to suppose that we know one point $h_q$ in the interior of
each $\Csetvar_q$ and we use \cororef{convexcone} in a transformed space where $h_q$ is mapped to the $z$-axis $(0, 1)$.
For this transformation we use the \emph{Householder reflection} \cite[Section~5.1.2]{golub2012matrix}
\[ \House{h} = I - \frac{2}{h^\Tr h} hh^\Tr. \]
The householder reflection is symmetric and orthogonal.

The optimization problem to solve is represented in \progref{quad}.
The transformation of this program to a semidefinite program can be done automatically using the using the standard Sum-of-Square procedure; see \cite{blekherman2012semidefinite}.

\begin{myprog}
  \label{prog:quad}
  \begin{align}
    \notag
    \max_{\substack{D_q \in \SymK, d_q \in \R^n,\\\delta_q \in \R, \lambda_{q \trs q'} \geq 0}}
    & \quad \sum_{q \in \Nodes} \log \det D_q
    \\
    %p_q(r(A_\sigma, c_\sigma)^\Tr x) & \geq \sum_{q \trs q' \in \Cs_\as} \lambda_{q \trs q'} p_{q'}(r(E_\sigma, 0)^{\Tr} x),\\
    \label{eq:quadl1}
    \begin{bmatrix}
      D_q & d_q\\
      d_q^\Tr & \delta_q + 1
    \end{bmatrix} & \succ 0\\
    \label{eq:quadp}
    p_q(y) & = y^\Tr \House{h_q}
    \begin{bmatrix}
      D_q & d_q\\
      d_q^\Tr & \delta_q
    \end{bmatrix}
    \House{h_q} y\\
    %p_q(r(A_\sigma, c_\sigma)^\Tr x) & \geq \sum_{q \trs q' \in \Cs_\as} \lambda_{q \trs q'} p_{q'}(r(E_\sigma, 0)^{\Tr} x),\\
    \label{eq:quads}
    p_q(r(A_\sigma, c_\sigma)^\Tr y) & \leq \lambda_{q \trs q'} p_{q'}(r(E_\sigma, 0)^{\Tr} y),\\
    \notag
    \forall q \in \Nodes, & \forall q \trs q', \forall y \in \R^{\nqx+1}\\
    \label{eq:quadc}
    p_q(-a, \beta) & \leq 0 \quad \forall q \in \Nodes, \forall a^\Tr x \leq \beta \text{ supporting } \Ps_q\\
    \label{eq:quadf}
    p_q(0, 1) & < 0 \quad \forall q \in \Nodes.
  \end{align}
\end{myprog}

The constraint \eqref{eq:quadl1} ensures both convexity of $\tau(\Csetvar_q)^*$ and the fact that $\det D_q$ does not overestimate the volume of the ellipsoid transformed by the Householder reflection.
The constraint \eqref{eq:quads} is the S-procedure applied to the condition \eqref{eq:dtahas}.
The constraint \eqref{eq:quadc} uses \eqref{eq:hs} to ensure that $\Csetvar_q$ is contained in $\Ps_q$.
The constraint \eqref{eq:quadf} ensures that $\tau(\Csetvar_q)^*$ has non-empty interior.
Note that if $\Ps_q$ has no unbounded subspace, \eqref{eq:quadf} is not necessary since the non-empty interior condition will already be ensured by \eqref{eq:quadc}.

\begin{mytheo}
  \label{theo:quad}
  Consider a \dtahas{} $\sys$ and points $(h_q \in \Ps_q)_{q \in \Nodes}$.
  The polynomial $p_q(x, z)$ is feasible for \progref{quad} if and only if
  there exists invariant convex sets $\Csetvar = (\Csetvar_q)_{q \in \Nodes}$
  such that $h_q \in \Csetvar_q$ for each $q \in \Nodes$
  and $\tau(\Csetvar_q)^*$ is the 0-sublevel set of $p_q(x, z)$.
  Moreover, the optimal solution of \progref{quad} is the solution that minimizes
  the sum of the logarithm of the volume of the intersection of the each cone $\tau(\Csetvar_q)^*$ with the hyperplane
  \( \{\, x \mid \langle h_q, x \rangle = 1 \,\}. \)
  \begin{proof}
    Consider a solution $p = (p_q(x, z))_{q \in \Nodes}$ of \progref{quad}.
    By \cororef{convexcone}, constraints \eqref{eq:quadl1} and \eqref{eq:quadp} are satisfied if and only if there exists ellipsoids $\Csetvar_q$ such that $\tau(\Csetvar_q)^*$ is the 0-sublevel set of $p_q(x, z)$.
    By \eqref{eq:hs}, constraint \eqref{eq:quadc} is satisfied if and only if $\Csetvar_q \subseteq \Ps_q$.
    By \lemref{sproc}, constraint \eqref{eq:quads} is satisfied if and only if
    \eqref{eq:dtahas} hold for all $q \trs q'$.
    Therefore, by \theoref{dtahas}, the solution $p$ is a feasible solution of \progref{quad} if and only if the sets $\Csetvar_q$ are invariant for $\sys$.

    Let $Q_q, c_q$ be such that $\Ellc{Q_q}{c_q} = \Ellq{D_q}{d_q}{\delta_q}$ and let $\lambda_q$ be such that $D_q = \lambda_q Q_q$.
    The volume of the intersection of $\tau(\Csetvar_q)^*$ with the hyperplane
    \( \{\, x \mid \langle h_q, x \rangle = 1 \,\} \)
    is $-\det(Q_q)$.
    Therefore, it remains to show that $\lambda_q = 1$ for an optimal solution.
    We observe that without the constraint~\eqref{eq:quadl1}, for any feasible solution, $D_q, d_q, \delta_q$ can be scaled by any positive constant
    while remaining feasible but affecting the objective function.
    By the Schur complement, constraint~\eqref{eq:quadl1} implies that
    \[ d_q^\Tr D_q^{-1} d_q - \delta_q \leq 1. \]
    Combining this inequality with equation~\eqref{eq:lambda} implies that $\lambda_q \leq 1$.
    Since the objective is to maximize $\det(D_q) = \lambda_q \det(Q_q)$,
    we know that if $(D_q, d_q, \delta_q)$ is optimal,
    then $\lambda_q = d_q^\Tr D_q^{-1} d_q - \delta_q = 1$.
  \end{proof}
\end{mytheo}

\begin{figure}[!h]
  \centering
  \includegraphics[width=0.4\textwidth]{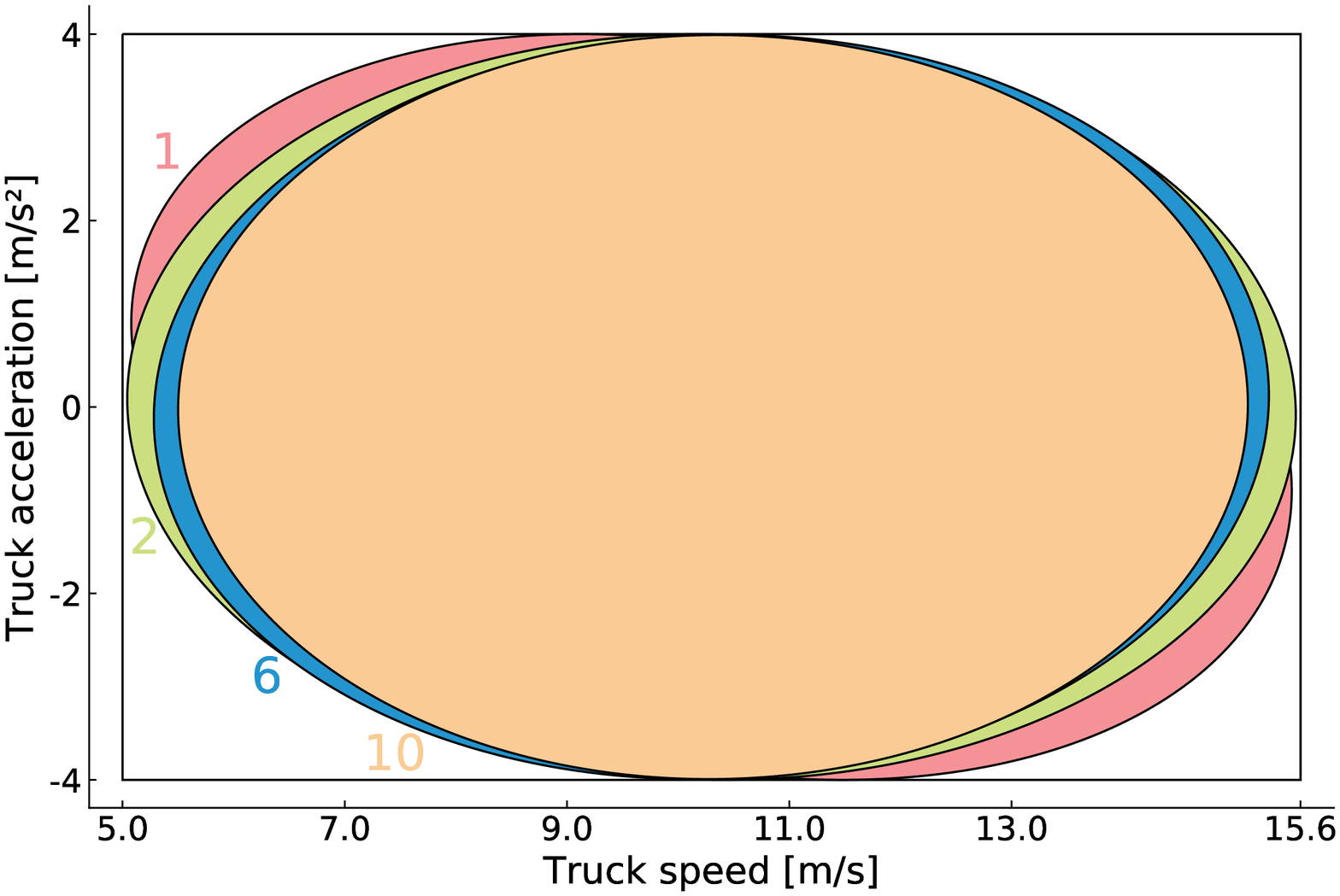}
  \caption{Projection onto the state $v_0$ and input $u$ of the optimal solution \new{of \progref{quad}} for \exemref{cruise3} at \node{} $q_{a0}$ for various numbers of trailers.}
  \label{fig:cruisesets}
\end{figure}

\begin{figure}[!h]
  \centering
  \includegraphics[width=0.4\textwidth]{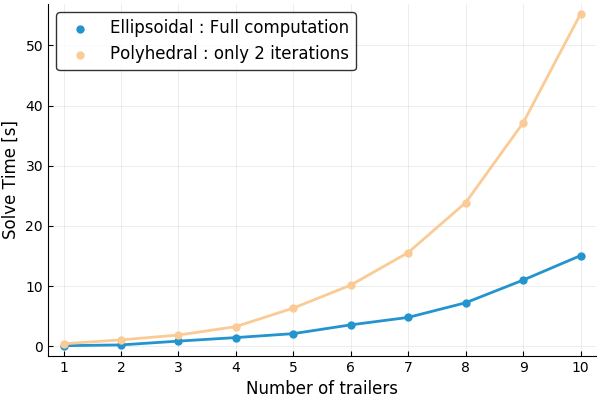}
  \caption{Computation time with Mosek 8.1.0.34 for \exemref{cruise3} \new{with various numbers of trailers compared
  to two iterations of the polyhedral approach (see e.g., the procedure p.~201 in \cite{blanchini2015set}) implemented with \new{the CDD library} \cite{fukuda1999cdd}.
  Note that after two iterations, the polyhedral sets obtained are not controlled invariant.
  One needs to wait for the convergence of the algorithm to obtain a controlled invariant set.
  Moreover, iterations are usually increasingly slower as the number of facets of the polyhedral sets increases with the iterations.}
  }
  \label{fig:cruisebench}
\end{figure}

\begin{myexem}
  \label{exem:cruise3}
  We apply \progref{quad} to \exemref{cruise2} with the same values for the parameters as the ones used in \cite{rungger2013specification},
  that is, $m_0 = \SI{500}{\kilogram}$, $m = \SI{1000}{\kilogram}$, $k_d = \SI{4600}{\newton\second\per\meter}$ and $k_s = \SI{4500}{\newton\per\kilogram}$.
  The values used for $h_q$ are the same for each \node{} $q \in \Nodes$: $u = d_i = 0$ and $v_0 = v_i = (5+v_a)/2$ for $i = 1, \ldots, M$. % $d$ and $u$, we use 0 and for the entries corresponding to the speed of the truck or a trailer we use $(5+v_a)/2$.
  %We observe that the optimal values of $\lambda_q$ is 1 for each $q \in \Nodes$.

  We vary the number of trailers $M$ from 1 to 10.
  \figref{cruisesets} represents the controlled invariant set at \node{} $q_{a0}$.
  As we can see, the constraints on the trailers are propagated to the truck and, as the number $M$ increases, the truck speed and acceleration become more constrained.

  The time taken by Mosek 8.1.0.34 (\cite{mosek2017mosek81034}) to solve the problem is given by \figref{cruisebench}\footnote{\new{We set $\lambda_{q \trs q'}$ to 1 for each transition $q \trs q'$ to make the problem convex}.}.
\end{myexem}

\section{Application to Model Predictive Control}
\label{sec:mpc}

As mentioned in the introduction, the controlled invariant sets can be used to derive a feedback control law.
We illustrate this with a Model Predictive Control (MPC) numerical experiment.
We consider a truck with one trailer ($M = 1$) as in \exemref{cruise3}.
The truck starts with speeds $v_0 = v_1 = \SI{10}{\meter\per\second}$ and spring displacement $d = \SI{0}{\meter}$
and has as objective to maximize the distance covered in \new{\SI{60}{\second}}.
The maximal speed is initially \SI{35}{\meter\per\second} but after \new{\SI{30}{\second}}, it drops to $v_a = \SI{15.6}{\meter\per\second}$.

In a classical MPC controller, the truck acceleration $u$ is controlled by solving a constrained optimal control problem up to horizon $H$.
%In the context of \exemref{cruise3}, the truck has \SI{0.8}{\second} to decrease its speed once it sees a new speed limitation.
%In this experiment, the MPC controller has access to the future safe sets up to time $H$ so it detects the speed limitation drop a time $H$ before it applies.
We observe that if $H \leq \new{\SI{9.2}{\second}}$, the controller is at some point unable to find values of $u$ satisfying input constraints
such that the state remains in the safe set.
%However, if $H \geq \new{\SI{9.6}{\second}}$, the controller is able to find input values that keep the state inside the safe set for \new{\SI{60}{\second}}.
%Note that the final state is not guaranteed to be a state from which the system can remain in the safe set after the \new{\SI{60}{\second}}.

\begin{figure}[!h]
  \centering
  \includegraphics[width=0.49\textwidth]{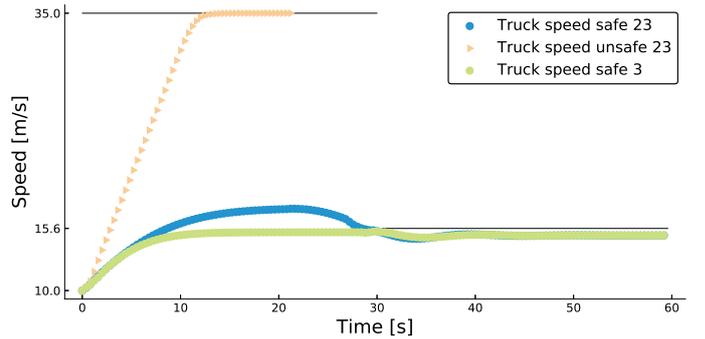}
  \caption{
    \new{Evolution with time of the speed of the truck for various MPC strategies.
    In the legend, \emph{safe} designates our MPC strategy using our computed invariant sets, while \emph{unsafe} designates a classical MPC approach.}
    %In the legend, \emph{optimal} refers to the optimal solution of the control problem and otherwise the number is the length of the horizon in time steps of \SI{0.4}{\second}.
    %The vertical line shows the switching time for the speed limitation and the horizontal lines show the current speed limitations.
    \new{The piecewise horizontal line represents the speed limitation at time $t$}.
    One can see that the MPC approach with invariant sets allows to remain in the safe set even with an horizon of 3 time steps.
    \new{Moreover, the unsafe controller can fail to find feasible values, as shown in \figref{acceleration}.}
  }
  \label{fig:speed}
\end{figure}

\begin{figure}[!h]
  \centering
  \includegraphics[width=0.49\textwidth]{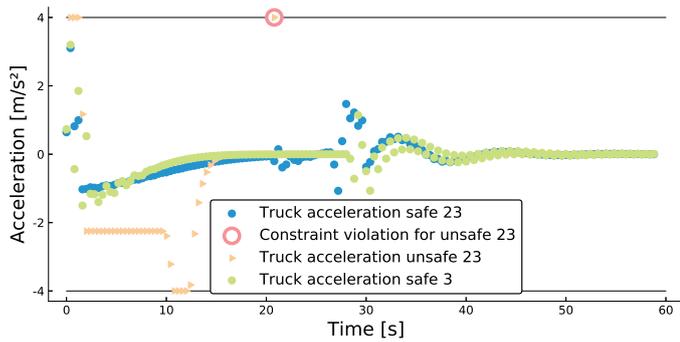}
  \caption{Acceleration of the truck in safe and unsafe mode. See \figref{speed} for the legend syntax.
    \new{We see (just after $t = \SI{20}{\second}$) that the unsafe controller requires a value $u > \SI{4}{\meter\per\second\squared}$ in order to remain in the safe set.}
    \new{Moreover}, we can see that the control is smoother in safe mode. 
    \new{Note that} using a longer horizon $H$ \new{renders the control even smoother, see e.g., between \SI{30}{\second} and \SI{40}{\second}}.
  }
  \label{fig:acceleration}
\end{figure}

For safety-critical applications, this lack of guarantee is not acceptable as it is necessary to be certain that the system can remain in the safe set.
Moreover, in a real-time context, the need to pick a large horizon is problematic as it increases the cost of online computations.
\new{In our setting, we constrain the state to remain in} the controlled invariant sets computed in \exemref{cruise3}\footnote{
  \new{\exemref{cruise3} corresponds to an MPC controller of horizon \SI{0.8}{\second}.
  An MPC controller of different horizon computes different controlled invariant sets by updating the hybrid system accordingly.}
} \new{and thereby} solve both issues.
Indeed, safety is guaranteed for arbitrarily long simulations and the length of the horizon does not influence safety so smaller length can be used.
Note that the controlled invariant sets can be computed offline so if it allows to reduce the horizon length,
it enables online computational cost to be moved offline.
Besides, constraining the state variables to belong to the ellipsoidal controlled invariant sets
is straightforward\footnote{The membership to $\Ellc{Q}{c}$ is second order cone representable. Indeed consider a Cholesky factorization $Q = L^\Tr L$, the inequality $(x-c)^\Tr Q (x-c) \leq 1$ can be rewritten as $\|L(x-c)\|_2 \leq 1$ where $\|\cdot\|_2$ is the Euclidean norm.}.
The results of the experiment can be found in \figref{speed} and \figref{acceleration}.
%\footnote{Note that the speed is lower during the 15 first seconds in the safe mode.
%  The reason is that when computing the controlled invariant sets,
%  the speed of the truck is constrained to be able to decrease in $\SI{0.8}{\second}$ to a change of speed limitation
%  while we allow the MPC controller to know that no speed limit change will occur during a time $H$ which enables
%  the truck to use faster speed in unsafe mode.}

%The goal of the experiment is to see how the result the controlled invariant sets computed in \exemref{cruise3}
%can be used and what is gained by using them.

\section{Conclusion}
We have developed a methodology for computing controlled invariant sets of Discrete-Time Affine Hybrid Control System (\dtahcs{}) and Discrete-Time Affine Hybrid Algebraic System (\dtahas{})
with \emph{autonomous switching} (see \remref{autonomous}).
This method can be \new{combined with semidefinite programming in order}
to compute ellipsoidal controlled invariant sets.
%The invariant sets obtained can be used in safety critical model predictive control applications
%to ensure safety while allowing to use smaller horizon hence decreasing online computational cost by precomputing controlled invariant sets.
We have shown that our technique can be used as a building block in a model predictive control scheme.
This allows, among other things, to reduce the online computational cost by precomputing controlled invariant sets.

We feel that we have only scratched the surface of the potential of the duality correspondence of \secref{duality}.
Many extensions of this work are possible
%such as hybrid systems with controlled switching and controlled invariant sets that are the sublevel sets of polynomial of higher degree.
such as hybrid systems with controlled switching, or the use of Sum-Of-Squares techniques in order to enrich the geometry of the possible invariant sets.

\new{The reformulation of the computation of controlled invariant sets of hybrid control system to the computation of invariant sets of hybrid algebraic system with \lemref{liftu} and \lemref{proju} allows to have a more behavioral invariance relation.
In the future, we would like to put our result in the framework of behavioral theory in order to investigate how to further generalize them; see \cite{willems2013introduction}.}

\bibliography{biblio}

\end{document}